\documentclass[11pt]{article}

\usepackage[margin=1in]{geometry}
\usepackage{amsmath,amsthm,amssymb}
\usepackage{commath,mathtools}
\usepackage{algorithm,algorithmic}
\usepackage{graphicx}
\usepackage{caption}
\usepackage{subcaption}
\usepackage{hyperref}
\usepackage{color}
\usepackage{rotating}
\usepackage{listings,booktabs}
\usepackage[authoryear,round,longnamesfirst]{natbib}
\usepackage[normalem]{ulem} %
\usepackage{bm}

\providecommand{\keywords}[1]{{\small \textbf{\textit{Keywords---}} #1}}

\DeclareMathOperator*{\argmin}{\mathrm{arg\,min}}

\title{A Machine Learning Enhanced Algorithm for the Optimal \\Landing Problem}

\author{
	Yaohua Zang\thanks{Zhejiang University, \em{yhchuang@zju.edu.cn}} \and Jihao Long\thanks{Princeton University, \em{jihaol@princeton.edu}} \and Xuanxi Zhang\thanks{Peking University, \em{zhangxuanxi@pku.edu.cn}} \and Wei Hu\thanks{Princeton University, \em{wh11@princeton.edu}} \and Weinan E\thanks{Peking University and Princeton University, \em{weinan@math.princeton.edu}} \and Jiequn Han\thanks{Flatiron Institute, \em{jiequnhan@gmail.com}}
}

\date{}

\begin{document}
\maketitle
\begin{abstract}
	We propose a machine learning enhanced algorithm for solving the optimal landing problem. Using Pontryagin's minimum principle, we derive a two-point boundary value problem for the landing problem. The proposed algorithm uses deep learning to predict the optimal landing time and a space-marching technique to provide good initial guesses for the boundary value problem solver. The performance of the proposed method is studied using the quadrotor example, a reasonably high dimensional and strongly nonlinear system. Drastic improvement in reliability and efficiency is observed.
\end{abstract}

\keywords{Optimal control, landing problem, deep neural networks, machine learning-based warm start}

\section{Introduction}
\label{sec:introduce}
The optimal landing problem is concerned with optimally controlling aerial vehicles to land on the target position. 
Developing effective numerical algorithms for the optimal landing problem has been a challenging task for some time, 
due to the high dimensionality of the state space and nonlinearity of the dynamics.
Traditionally, there are two approaches for solving optimal control and landing problems: the direct method and the indirect method.
The direct method \citep{vanderbei1999loqo,ross2002direct} first translates the optimal control problem to a nonlinear optimization problem by discretizing time and then solves the discretized problem using well-developed optimization solvers. This approach has been applied to lunar landing~\citep{liu2006nonlinear}, rocket landing~\citep{liu2019fuel}, quadrotor landing~\citep{bouktir2008trajectory,hu2017time}, to name a few.
The indirect method \citep{bock1984multiple,bertolazzi2005symbolic,wang2009solving} is based on the Pontryagin Minimum/Maximum Principle (PMP), which can be interpreted as the first-order optimality conditions of the optimal control problem. 
One then solves the two-point boundary value problem (TPBVP) derived from the PMP.
The indirect method has also been applied to a variety of landing problems with different terminal constraints~\citep{guo2011optimal,assellaou2016hamilton,hu2015fast}.
Despite the progress, these methods still have serious limitations such as %
the sub-optimality of the trajectory \citep{foehn2021time,romero2021model}, reliance on good initial guesses \citep{geisert2016trajectory}, and long computation time.

In recent years, starting from \cite{han2016deep}, deep neural networks (DNN) have been widely used to solve high dimensional optimal control 
problem  (see, e.g., \cite{tang2018learning,tang2019data,zhu2019deep,beppu2021value,nakamura2021adaptive}).
\cite{sanchez2018real} used the DNNs to approximate the mapping from states to optimal actions and presented applications on several landing problems.
\cite{tang2018learning} introduced a trajectory optimization approach that achieved real-time performance by combining machine learning to predict optimal trajectories with refinement by quadratic optimization.
\cite{zhu2019deep} used DNNs to learn the optimal action to improve the computational efficiency for the fuel-optimum lunar landing problem.
\cite{Shi2019Neural} presented a deep learning-based robust nonlinear controller to improve the control performance of a quadrotor during landing.
\cite{you2020learning} developed a learning-based optimal control method for the Mars entry and power descent guidance to find optimal guidance laws in real-time.
Instead of learning the feedback control directly,
\cite{nakamura2021adaptive} proposed to learn the value function via DNNs and predict the optimal feedback through the dynamic programming principle.
Although they have shown great potential for the optimal control of high dimensional and strongly nonlinear systems, current DNN-based algorithms are still not robust enough~\citep{chen2018optimal,nakamura2021neural}.
One major obstacle is the sensitive dependence on a good initialization at the various stages of the algorithms.

This paper proposes a new numerical method that combines the traditional indirect method and DNNs to solve the optimal landing problems with much improved reliability and efficiency. Overall we will use DNN-based prediction of the terminal time and space-marching technique to warm start the solving process and accelerate convergence.
We will take the quadrotor unmanned aerial vehicles (UAVs) as an example to demonstrate the methodology. Rotary-wing UAVs have received widespread attention in recent years due to their wide range of application scenarios, including package delivery, film photography, agricultural inspections, and search and rescue missions.
Among the many types of UAVs, the quadrotor is the most prominent aerial system. 
Due to its simplicity and versatility, the quadrotor has become the most flexible and maneuverable drone \cite{ackermann2020ai,verbeke2018experimental}.
We consider the full quadrotor dynamic model and aim to achieve an optimal landing path with minimum time and control effort.
We start from the indirect method, utilizing the PMP to transform the original optimal landing problem into a TPBVP.
One critical issue of the TPBVP solver is to find good initial guesses \citep{tsiotras2011initial,nakamura2021adaptive}.
To overcome this difficulty, we design a DNN-based algorithm to provide an initial guess of the optimal landing time and a space-marching scheme to provide an initial guess of the solution.
Compared to the baseline TPBVP solver, the proposed algorithm obtains the optimal landing trajectory with a much higher success rate and less computation time.

The paper is organized as follows. 
Section \ref{sec:problem} presents the general form of control problems that our algorithm will be applied to and the specification of the optimal landing problem for the quadrotor. 
Our machine learning enhanced algorithm will be introduced in Section \ref{sec:DNN}, with numerical results demonstrating its efficiency. 
Section \ref{sec:MPC} presents the difference between the obtained optimal solution and the suboptimal solution provided by the model predictive control (MPC) method, another popular method for solving optimal control problems.
Section \ref{sec:summary} concludes the paper with some discussions on future work.

\section{Formulation of the optimal control problem}
\label{sec:problem}
We consider a deterministic system defined by the following ordinary differential equation (ODE)
\begin{equation}\label{eq:dynamic}
  \begin{cases}
    \dot{\bm{x}}(t) = f(\bm{x}(t), \bm{u}(t)), \quad t\in[0,t_f] \\
     \bm{x}(0) = \bm{x}_0, \; \bm{g}(\bm{x}(t_f)) = \bm{0},
  \end{cases}
\end{equation}
where $\bm{x}(t)\in \mathbb{R}^n$ indicates the states, $\bm{u}(t)\in \mathcal{U}\subset \mathbb{R}^m$ represents the controls with $\mathcal{U}$ being the admissible set of the controls, $\bm{f}:\mathbb{R}^n\times\mathcal{U}\mapsto \mathbb{R}^n$ and $\bm{g}: \mathbb{R}^n \mapsto \mathbb{R}^k$ are smooth functions describing the dynamics and terminal condition. We call $\{\bm{x},\bm{u},t_f\}$ a feasible path if \eqref{eq:dynamic} is satisfied and use $\mathcal{P}$ to denote the set of all feasible paths.
The total cost is defined as
\begin{equation}\label{eq:cost}
  \mathcal{C}[\bm{x},\bm{u},t_f]= \int^{t_f}_{0} L(\bm{x}(t),\bm{u}(t)) dt,
\end{equation}
where $L: \mathbb{R}^n\times\mathcal{U} \mapsto \mathbb{R}$ is the running cost, which is assumed to be smooth and non-negative.

We will consider two different but closely related problems. In the first problem, $t_f$ is a given positive constant and we aim to minimize the performance function over all feasible paths with a fixed terminal time $t_f$:
\begin{equation}\label{fix_pro}
    \min_{(\bm{x},\bm{u}):\{\bm{x},\bm{u},t_f\} \in \mathcal{P}} \mathcal{C}[\bm{x},\bm{u},t_f].
\end{equation}
We call this problem a \textit{fixed terminal time problem}. 
We are also interested in the  \textit{free terminal time problem}, where we aim to minimize the performance function over all feasible paths:
\begin{equation}\label{free_pro}
    \min_{\{\bm{x},\bm{u},t_f\} \in \mathcal{P}} \mathcal{C}[\bm{x},\bm{u},t_f].
\end{equation}
Our ultimate goal is to solve this free terminal time problem to find the optimal landing trajectory and corresponding control with minimum time and control effort.

We remark that, in more complicated cases, there can be path constraints such as $c(\bm{x}(t), \bm{u}(t))\geq 0$ in the dynamical system.
Here we only consider problems free of path constraints to highlight the main features of the techniques we introduce. 
These techniques can be naturally extended to problems with path constraints  and, we intend to study this in future work.

\subsection{Pontryagin's Minimum Principle}
Pontryagin's Minimum Principle (PMP) establishes a set of  necessary conditions for optimality, which converts the optimal control problem \eqref{fix_pro} or \eqref{free_pro} to two-point boundary value problems (TPBVP). 
Algorithms based on the PMP are usually called  indirect methods~\cite{bock1984multiple,bertolazzi2005symbolic,wang2009solving}.
To simplify the discussion, we assume that both the fixed terminal time problem \eqref{fix_pro} and the free terminal time problem \eqref{free_pro} admit a unique minimizer and the solutions of corresponding TPBVPs are unique. This assumption ensures the optimality of the solutions of TPBVPs. In the  examples considered here, we have not found multiple solutions to the TPBVPs.

To state the PMP, we introduce a costate variable $\bm{\lambda} \in \mathbb{R}^n$ and define the Hamiltonian
\begin{equation}\label{eq:H}
  H(\bm{x}, \bm{\lambda}, \bm{u})= L(\bm{x}, \bm{u})+\bm{\lambda} \cdot f(\bm{x},\bm{u}).
\end{equation}
The PMP reduces the fixed terminal time problem \eqref{fix_pro} to a system of equations in the form of
\begin{equation}\label{eq:bvp}
  \begin{cases}
    \dot{\bm{x}}(t)= \partial^T_{\bm{\lambda}} H(\bm{x}(t), \bm{\lambda}(t), \bm{u}^*(t)),\\
    \dot{\bm{\lambda}}(t)= -\partial^T_{\bm{x}} H(\bm{x}(t), \bm{\lambda}(t), \bm{u}^*(t)),
  \end{cases}
\end{equation}
together with the boundary conditions given by the original one augmented with the transversality conditions:
\begin{equation}\label{eq:bvp_bd}
  \begin{cases}
    \bm{x}(0) = \bm{x}_0, \\
  \bm{g}(\bm{x}(t_f)) = \bm{0}, \\
  \nabla \bm{g}(\bm{x}(t_f))\bm{\alpha} =  \bm{\lambda}(t_f).
  \end{cases}
\end{equation}
Here $\nabla \bm{g}(\bm{x}(t_f)) \in \mathbb{R}^{n\times k}$, and $\bm{\alpha} \in \mathbb{R}^k$ is a multiplier.
In addition, the optimal control $\bm{u}^*(t)$ should satisfy the minimization of the Hamiltonian at each $t$:
\begin{equation}\label{eq:optimal_u}
  \bm{u}^*(t) = \argmin_{\bm{u}\in \mathcal{U}} H(\bm{x}, \bm{\lambda}, \bm{u}).
\end{equation}
Equations \eqref{eq:bvp}, \eqref{eq:bvp_bd} and \eqref{eq:optimal_u} together complete the PMP for the fixed terminal time problem \eqref{fix_pro}. 
For the free terminal time problem \eqref{free_pro}, besides \eqref{eq:bvp}, \eqref{eq:bvp_bd} and \eqref{eq:optimal_u},
we need the extra condition for the optimal terminal time $t_f$:
\begin{equation}
    H(\bm{x}(t_f),\bm{\lambda}(t_f),\bm{u}^*(t_f)) = 0.
\end{equation}
See \cite{hartl1995survey} for the proof of the above PMP and detailed discussions.

\subsection{The optimal landing problem}
\subsubsection{The full dynamic model of quadrotor}
To introduce the dynamics of the quadrotor, we let $\{\mathcal{O}_E, \mathcal{X}_E, \mathcal{Y}_E, \mathcal{Z}_E\}$ denote the Earth-fixed coordinate system and $\{\mathcal{O}_b, \mathcal{X}_b, \mathcal{Y}_b, \mathcal{Z}_b\}$ the body-fixed coordinate system, whose origin $\mathcal{O}_b$ is at the center of mass (CoM) of the quadrotor.
Then the dynamics of the quadrotor can be modeled as follows
\begin{equation}\label{eq:uav_old}
  \begin{cases}
    \dot{\bm{p}}= \bm{R}^T(\bm{\eta})\bm{v}_b\\
    \dot{\bm{v}}_b= - \bm{w}_b\times \bm{v}_b - \bm{R}(\bm{\eta})\bm{g} + \frac{1}{m}\bm{f}_u\\
    \dot{\bm{\eta}} = \bm{K}(\bm{\eta})\bm{w}_b \\
    \dot{\bm{w}}_b = - \bm{J}^{-1}\bm{w}_b\times \bm{J}\bm{w}_b + \bm{J}^{-1}\bm{\tau}_u,
  \end{cases}
\end{equation}
where $\bm{p}= (x, y, z)^T$ is the inertial position of the CoM in the Earth-fixed coordinates and $\bm{v}_b=(v_x, v_y, v_z)$ is the linear velocity of the quadrotor expressed in the body-fixed coordinates. $\bm{\eta} = (\phi, \theta, \psi)$ is the attitude of the quadrotor in the Earth-fixed coordinates defined by the Euler angles: roll($\phi$), pitch($\theta$) and yaw($\psi$). $\bm{w}_b = (p, q, r)^T$ denotes the angular velocity in the body-fixed coordinates.
In total $\bm{x} = (\bm{p}^T, \bm{v}_b^T, \bm{\eta}^T, \bm{w}_b^T)^T \in \mathbb{R}^{12}$ denotes the state variable. 
$\bm{f}_u= (0, 0, T)^T$ and $\bm{\tau}_u = (\tau_x, \tau_y, \tau_z)^T$ are the total thrust and body torques from the four rotors, which are forces applied by the control variables to adjust the quadrotor's dynamics. The constants $m$ and $\bm{g}= (0,0, g)^T$ denote the mass and the gravity vector ($g=9.81 m/s^2$ denotes the acceleration of gravity on Earth), respectively.
The direction cosine matrix $\bm{R}(\bm{\eta})\in SO(3)$ (representing the transformation from the Earth-fixed coordinates to the body-fixed coordinates), attitude kinematic matrix $\bm{K}(\bm{\eta})$ (relating the time derivative of the attitude representation with the associated angular rate) and constant inertia matrix $\bm{J}$ are defined as follows
\begin{equation*}
  \bm{R}(\bm{\eta}) =
  \begin{bmatrix}
    \cos{\theta}\cos{\psi} & \cos{\theta}\sin{\psi} & -\sin{\theta}\\
    \sin{\theta}\cos{\psi}\sin{\phi}-\sin{\psi}\cos{\phi} &
    \sin{\theta}\sin{\psi}\sin{\phi}+\cos{\psi}\cos{\phi} &
    \cos{\theta}\sin{\phi} \\
    \sin{\theta}\cos{\psi}\cos{\phi}+\sin{\psi}\sin{\phi} &
    \sin{\theta}\sin{\psi}\cos{\phi}-\cos{\psi}\sin{\phi} &
    \cos{\theta}\cos{\phi}
  \end{bmatrix},
\end{equation*}

\begin{equation*}
  \bm{K}(\bm{\eta}) =
  \begin{bmatrix}
    1 & \sin{\phi}\tan{\theta} & \cos{\phi}\tan{\theta}\\
    0 & \cos{\phi} & -\sin{\phi} \\
    0 & \sin{\phi}\sec{\theta} & \cos{\phi}\sec{\theta}
  \end{bmatrix},
\end{equation*}

\begin{equation*}
  \bm{J} = \text{diag}(J_{x}, J_{y}, J_{z}),
\end{equation*}
where $J_{x}, J_{y}$, and $J_{z}$ are the moments of inertia of the quadrotor in the $x$-axis, $y$-axis, and $z$-axis, respectively.

To ease the notation, we denote the controls as $\bm{u}=(T, \tau_x, \tau_y, \tau_z)^T$. Then we have $\bm{f}_u = A\bm{u}$ and $\bm{\tau}_u = B\bm{u}$ with $A$ and $B$ defined as
\begin{equation*}
  A
  =
  \begin{bmatrix}
    0 & 0 & 0 & 0 \\
    0 & 0 & 0 & 0 \\
    1 & 0 & 0 & 0
  \end{bmatrix},
  \quad\quad
  B
  =
  \begin{bmatrix}
    0 & 1 & 0 & 0 \\
    0 & 0 & 1 & 0 \\
    0 & 0 & 0 & 1
  \end{bmatrix}.
\end{equation*}
Note that in practice the quadrotor is directly controlled by the the individual rotor thrusts $\bm{F}=(F_1, F_2, F_3, F_4)^T$, and we have the relation $\bm{u}=E\bm{F}$ with
\begin{equation*}
  E
  =
  \begin{bmatrix}
    1 & 1 & 1 & 1  \\
    0 & l & 0 & -l \\
    -l & 0 & l & 0 \\
    c & -c & c & -c
  \end{bmatrix},
\end{equation*}
where $l$ is the distance from the rotor to the UAV's center of gravity and $c$ is a constant that relates the rotor angular momentum to the rotor thrust (normal force). 
So once we obtain the optimal control $\bm{u}^*$, we are able to get the optimal $\bm{F}^*$ immediately by the relation $\bm{F}^*=E^{-1}\bm{u}^*$. We also introduce the skew-symmetric matrix $S(\bm{w}_b)$
\begin{equation*}
  \bm{S}(\bm{w}_b) =
  \begin{bmatrix}
    0 & r & -q\\
    -r & 0 & p \\
    q & -p & 0
  \end{bmatrix},
\end{equation*}
which has the property
\begin{equation*}
  -\bm{w}_b\times \bm{v}_b = \bm{S}(\bm{w}_b)\bm{v}_b \quad \text{and}\quad - \bm{J}^{-1}\bm{w}_b\times \bm{J}\bm{w}_b = \bm{J}^{-1}\bm{S}(\bm{w}_b)\bm{J}\bm{w}_b.
\end{equation*}
Then equation \eqref{eq:uav_old} can be finally rewritten as
\begin{equation}\label{eq:uav}
  \begin{cases}
    \dot{\bm{p}}= \bm{R}^T(\bm{\eta})\bm{v}_b\\
    \dot{\bm{v}}_b= \bm{S}(\bm{w}_b)\bm{v}_b - \bm{R}(\bm{\eta})\bm{g} + \frac{1}{m}A\bm{u}\\
    \dot{\bm{\eta}} = \bm{K}(\bm{\eta})\bm{w}_b \\
    \dot{\bm{w}}_b = \bm{J}^{-1}\bm{S}(\bm{w}_b)\bm{J}\bm{w}_b + \bm{J}^{-1}B\bm{u}.
  \end{cases}
\end{equation}

\subsubsection{The optimal control problem}
We aim to solve the landing problem with minimum control effort and shortest landing time under the dynamics described in \eqref{eq:uav}.
That is, to find the optimal controls to steer the quadrotor from some initial states $\bm{x}_0\in \mathcal{S}_0$ to a target state $\bm{x}_{t_f}\in \mathcal{S}_{T}\coloneqq \{\bm{x}: \bm{g}(\bm{x})=\bm{0}\}$. 

For the landing problem, the terminal set $\mathcal{S}_{T}$ has the form (recalling $\bm{x} = (\bm{p}^T, \bm{v}_b^T, \bm{\eta}^T, \bm{w}_b^T)^T$)
\begin{equation*}
\{\bm{x}_{t_f}\;|\; \bm{p}(t_f)=\bm{v}(t_f)=\bm{w}(t_f) =\bm{0},\, \phi(t_f)=\theta(t_f)=0\}.
\end{equation*}
The running cost $L$ in \eqref{eq:cost} is given by
\begin{align*}
  L(\bm{x}, \bm{u}) = 1 + (\bm{u}-\bm{u}_d)^TQ_u(\bm{u}-\bm{u}_d),
\end{align*}
where $\bm{u}_d=(mg,0,0,0)$ represents the reference control that balances with gravity and $Q_{u}=\text{diag}(1, 1, 1, 1)$ represents the weight matrix characterizing the cost of deviating from the reference control.
Then, the optimal control problem with fixed terminal time can be written as
\begin{equation}\label{eq:opt_problem}
  \begin{array}{ll}
    \min_{\bm{x},\bm{u}} & \mathcal{C}[\bm{x},\bm{u},t_f] \\ \\
    \mathop{subject}\ \  \mathop{to} &
    \begin{cases}
      \dot{\bm{p}}= \bm{R}^T(\bm{\eta})\bm{v}_b\\
      \dot{\bm{v}}_b= \bm{S}(\bm{w}_b)\bm{v}_b - \bm{R}(\bm{\eta})\bm{g} + \frac{1}{m}A\bm{u}\\
      \dot{\bm{\eta}} = \bm{K}(\bm{\eta})\bm{w}_b \\
      \dot{\bm{w}}_b = \bm{J}^{-1}\bm{S}(\bm{w}_b)\bm{J}\bm{w}_b + \bm{J}^{-1}B\bm{u} \\
      \bm{x}(0) = \bm{x}_0 \\
      \bm{p}(t_f) = \bm{v}(t_f) =  \bm{w}(t_f) = \bm{0},\,
      \phi(t_f) = \theta(t_f) = 0.
    \end{cases}
  \end{array}
\end{equation}
The optimal control problem with free terminal time needs to further find the optimal terminal time $t_f$ to minimize the total cost in \eqref{eq:opt_problem}.

To derive the TPBVP for the problem \eqref{eq:opt_problem}, we denote $\bm{\lambda} = (\bm{\lambda_p}^T, \bm{\lambda_v}^T, \bm{\lambda_{\eta}}^T, \bm{\lambda_{w}}^T)^T \in \mathbb{R}^{12}$ as the costate of $\bm{x}$ and write down the Hamiltonian
\begin{align*}
  H(\bm{x},\bm{\lambda}, \bm{u}) = &~ L(\bm{x}, \bm{u})+\bm{\lambda}_p \cdot (\bm{R}^T(\bm{\eta})\bm{v}_b) + \bm{\lambda}_v \cdot (\bm{S}(\bm{w}_b)\bm{v}_b - \bm{R}(\bm{\eta})\bm{g} + \frac{1}{m}A\bm{u})
  \\
  &~  + \bm{\lambda}_{\eta}\cdot (\bm{K}(\bm{\eta})\bm{w}_b) + \bm{\lambda}_w\cdot (\bm{J}^{-1}\bm{S}(\bm{w}_b)\bm{J}\bm{w}_b+\bm{J}^{-1}B\bm{u}).
\end{align*}
Therefore, the TPBVP for the optimal control problem \eqref{eq:opt_problem} is 
\begin{equation}\label{eq:bvp_landing}
  \begin{cases}
    \dot{\bm{p}}= \bm{R}^T(\bm{\eta})\bm{v}_b\\
    \dot{\bm{v}}_b= \bm{S}(\bm{w}_b)\bm{v}_b - \bm{R}(\bm{\eta})\bm{g} + \frac{1}{m}A\bm{u}^*\\
    \dot{\bm{\eta}} = \bm{K}(\bm{\eta})\bm{w}_b \\
    \dot{\bm{w}}_b = \bm{J}^{-1}\bm{S}(\bm{w}_b)\bm{J}\bm{w}_b + \bm{J}^{-1}B\bm{u}^* \\
    \dot{\bm{\lambda}}_p = \bm{0}\\
    \dot{\bm{\lambda}}_v = - \bm{R}(\bm{\eta})\bm{\lambda}_{p} - \bm{S}^T(\bm{w}_b)\bm{\lambda}_v\\
    \dot{\bm{\lambda}}_{\eta} = 
    -\frac{\partial(\lambda_p\cdot\bm{R}^T(\bm{\eta})\bm{v}_b)}{\partial\bm{\eta}}
    +\frac{\partial(\bm{\lambda}_{v}\cdot\bm{R}(\bm{\eta})\bm{g})}{\partial\bm{\eta}} 
    -\frac{\partial(\bm{\lambda}_{\eta}\cdot \bm{K}(\bm{\eta})\bm{w}_b)}{\partial\bm{\eta}}\\
    \dot{\bm{\lambda}}_w = 
    -\bm{K}^T(\bm{\eta})\bm{\lambda}_{\eta}
    -\frac{\partial(\bm{\lambda}_v\cdot \bm{S}(\bm{w}_b)\bm{v}_b)}{\partial\bm{w}_b}
    -\frac{\partial(\bm{\lambda}_w\cdot \bm{J}^{-1}\bm{S}(\bm{w}_b)\bm{J}\bm{w}_b)}{\partial\bm{w}_b},
  \end{cases}
\end{equation}
with boundary conditions
\begin{equation}\label{eq:bvp_bd_landing}
  \begin{cases}
    \bm{x}(0) = \bm{x}_0, \\
    \bm{p}(t_f)=\bm{v}(t_f)=\bm{w}(t_f)=\bm{0},\\ 
    \phi(t_f)=\theta(t_f)=\lambda_{\psi}(t_f)=0.
  \end{cases}
\end{equation}
For the optimal control problem with free terminal time, the optimal final time $t_f$ is further determined by
\begin{equation}\label{eq:bvp_tf_landing}
  \\
  H(\bm{x}(t_f), \bm{\lambda}(t_f), \bm{u}^*(t_f)) = 0.
\end{equation}
In both cases, the optimal feedback control $\bm{u}^*$ at state $\bm{x}$ and time $t$ is
\begin{align*}
  \bm{u}^* = &~\argmin_{\bm{u}}H(\bm{x}, \bm{\lambda},\bm{u}) \\
  = &~\bm{u}_d - (Q^T_u+Q_u)^{-1}(\frac{1}{m}A^T\bm{\lambda}_v+B^T(\bm{J}^T)^{-1}\bm{\lambda}_w).
\end{align*}
In the free terminal time problem, the undetermined $t_f$  brings difficulty to the time integrator in the TPBVP solver since the time interval keeps changing during the solving process. 
To circumvent this difficulty, we consider $t_f$ as a new state $x_{n+1}=t_f$ and introduce a new variable $\tau = t/t_f= t/x_{n+1}$.
We then have
\begin{equation}
  \dot{x}_{n+1}=0
\end{equation}
and
\begin{equation}\label{eq:bvp_free_a}
  \begin{cases}
    \dot{\bm{p}}= x_{n+1}(\bm{R}^T(\bm{\eta})\bm{v}_b)\\
    \dot{\bm{v}}_b= x_{n+1}(\bm{S}(\bm{w}_b)\bm{v}_b - \bm{R}(\bm{\eta})\bm{g} + \frac{1}{m}A\bm{u}^*)\\
    \dot{\bm{\eta}} = x_{n+1}(\bm{K}(\bm{\eta})\bm{w}_b) \\
    \dot{\bm{w}}_b = x_{n+1}(\bm{J}^{-1}\bm{S}(\bm{w}_b)\bm{J}\bm{w}_b + \bm{J}^{-1}B\bm{u}^*).
  \end{cases}
\end{equation}
The cost function becomes
\begin{equation*}
  \mathcal{C} = x_{n+1}\int^{1}_{0} L(\bm{x}(\tau), \bm{u}(\tau)) d\tau,
\end{equation*}
and the Hamiltonian becomes
\begin{align*}
  H(\bm{x}, x_{n+1}, \bm{\lambda}, \bm{u}) = &~  x_{n+1}\bigg(L(\bm{x}, \bm{u}) + \bm{\lambda}_p \cdot (\bm{R}^T(\bm{\eta})\bm{v}_b) + \bm{\lambda}_v \cdot (\bm{S}(\bm{w}_b)\bm{v}_b - \bm{R}(\bm{\eta})\bm{g} + \frac{1}{m}A\bm{u})\\
  &~ + \bm{\lambda}_{\eta}\cdot (\bm{K}(\bm{\eta})\bm{w}_b)  + \bm{\lambda}_w\cdot (\bm{J}^{-1}\bm{S}(\bm{w}_b)\bm{J}\bm{w}_b+\bm{J}^{-1}B\bm{u})\bigg).
\end{align*}
Similarly, the costate $\bm{\lambda}$ satisfies the following equation
\begin{equation}\label{eq:bvp_free_b}
  \begin{cases}
    \dot{\bm{\lambda}}_p = \bm{0}\\
    \dot{\bm{\lambda}}_v =  -x_{n+1} (\bm{R}(\bm{\eta})\bm{\lambda}_{p} + \bm{S}^T(\bm{w}_b)\bm{\lambda}_v)\\
    \dot{\bm{\lambda}}_{\eta} = - x_{n+1}(Q_{\eta}^T+Q_{\eta})(\bm{\eta}-\bm{\eta}_d)
    -x_{n+1}(\frac{\partial(\lambda_p\cdot\bm{R}^T(\bm{\eta})\bm{v}_b)}{\partial\bm{\eta}}
    -\frac{\partial(\bm{\lambda}_{v}\cdot\bm{R}(\bm{\eta})\bm{g})}{\partial\bm{\eta}} +\frac{\partial(\bm{\lambda}_{\eta}\cdot \bm{K}(\bm{\eta})\bm{w}_b)}{\partial\bm{\eta}})\\
    \dot{\bm{\lambda}}_w =
    -x_{n+1}(\frac{\partial(\bm{\lambda}_v\cdot \bm{S}(\bm{w}_b)\bm{v}_b)}{\partial\bm{w}_b}
    +\bm{K}^T(\bm{\eta})\bm{\lambda}_{\eta}
    +\frac{\partial(\bm{\lambda}_w\cdot \bm{J}^{-1}\bm{S}(\bm{w}_b)\bm{J}\bm{w}_b)}{\partial\bm{w}_b}).
  \end{cases}
\end{equation}
The transversality condition \eqref{eq:bvp_tf_landing} then becomes
\begin{equation}
\label{eq:bvp_free_c}
  H(\bm{x}(1), \bm{\lambda}(1), \bm{u}^*(1)) = 0.
\end{equation}
Therefore, we transform the TPBVP \eqref{eq:bvp_landing}-\eqref{eq:bvp_tf_landing} with undetermined time-horizon to a TPBVP with a fixed time-horizon $[0,1]$, which is easier to solve for a TPBVP solver.

\section{Machine learning enhanced algorithm}
\label{sec:DNN}
\begin{figure}[!htb]
	\centering
	\begin{subfigure}[b]{.45\textwidth}
	\includegraphics[width=\textwidth]{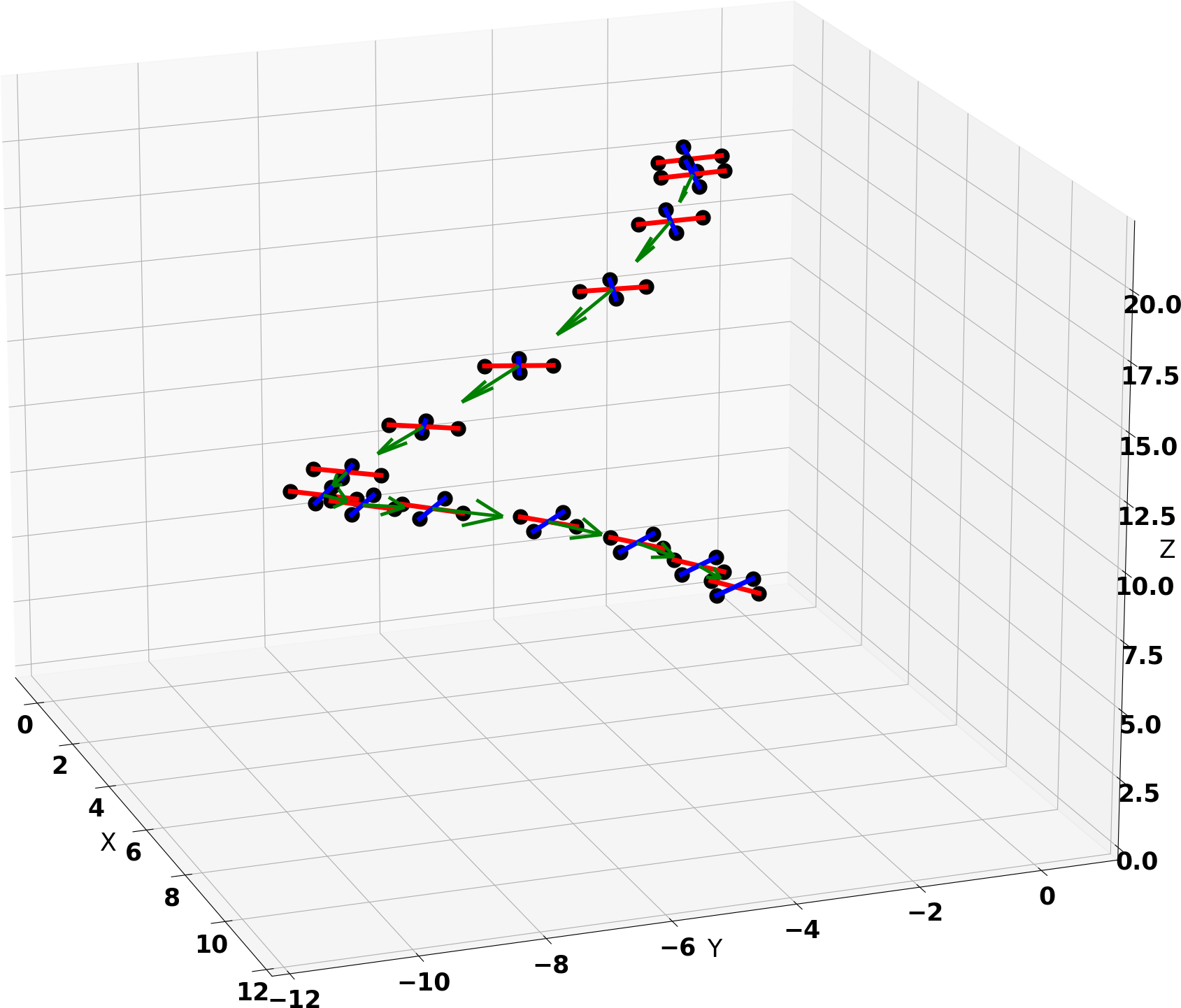}
	\vspace{-0.4em}
	\caption{}
	\label{subfig:path_3d}
	\end{subfigure}\\
	\begin{subfigure}[b]{.44\textwidth}
	\includegraphics[width=\textwidth]{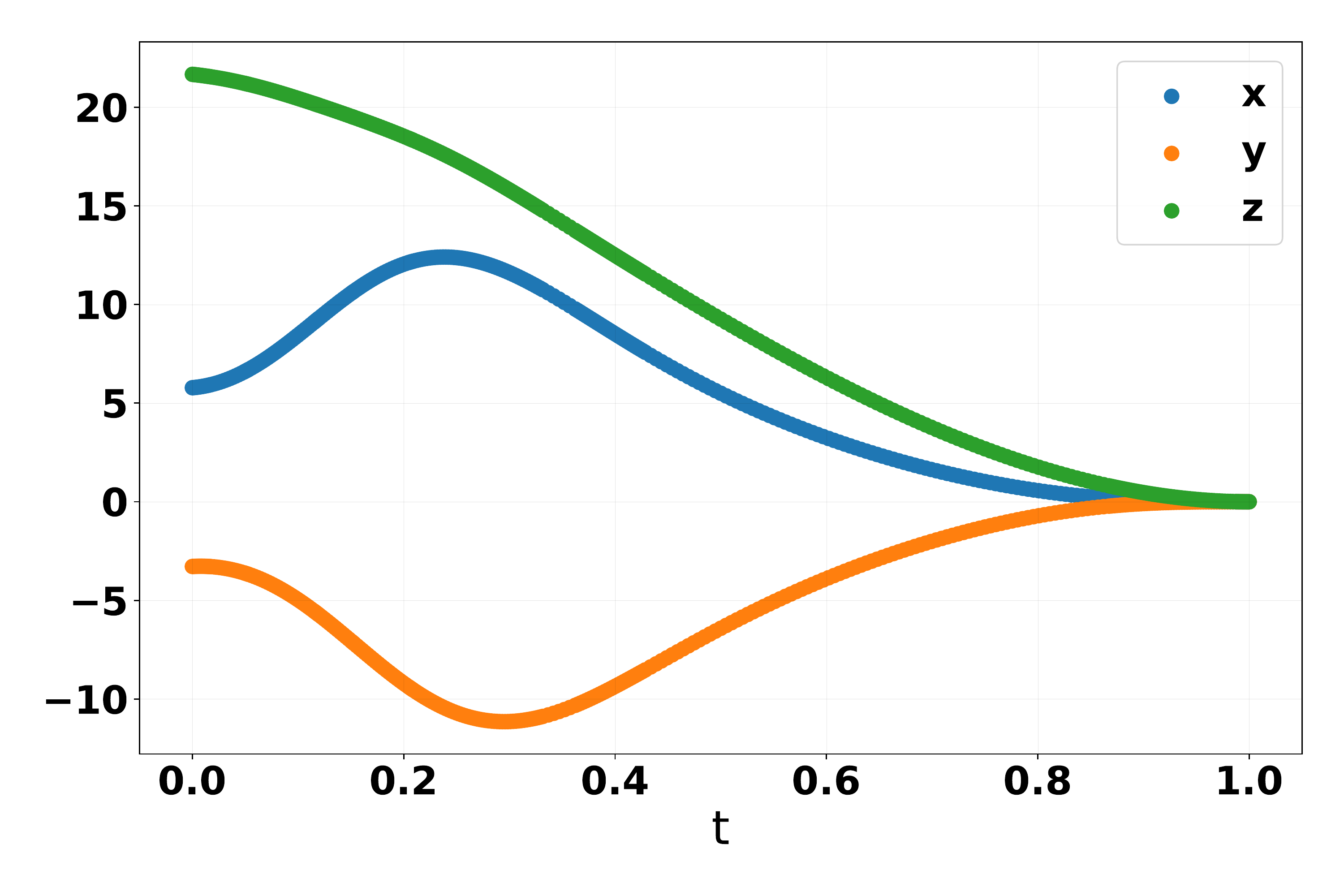}
	\caption{}
	\label{subfig:p_3d}
	\end{subfigure}
	\begin{subfigure}[b]{.44\textwidth}
	\includegraphics[width=\textwidth]{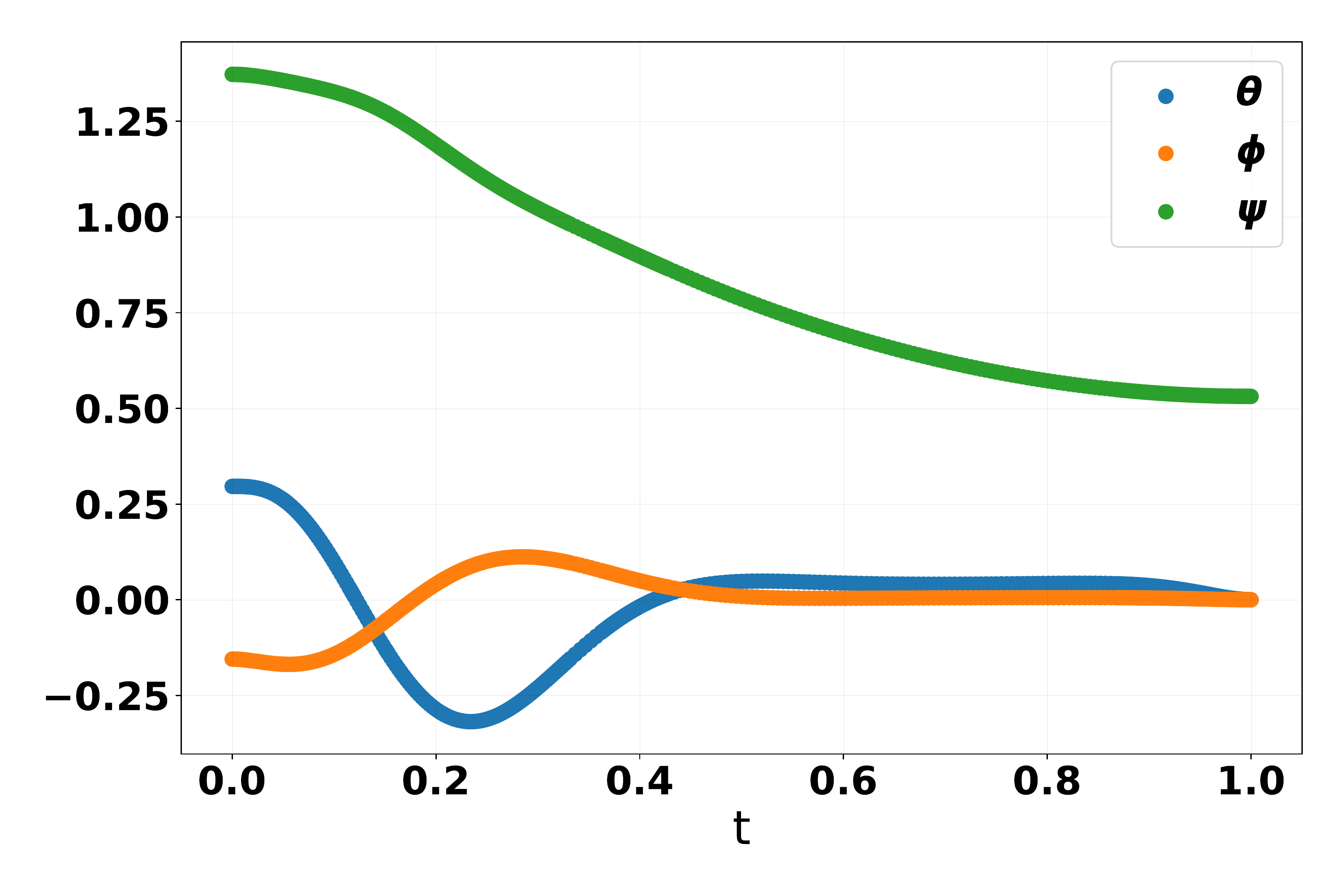}
	\caption{}
	\label{subfig:eta_3d}
	\end{subfigure}
	\caption{(a) An example of the optimal trajectory for the landing problem of quadrotor; (b) The obtained optimal $\bm{p}$ vs. time; (c) The obtained optimal $\bm{\eta}$ vs. time.}
	\label{fig:uav_3d}
\end{figure}

In this section, we present our algorithm for solving the optimal landing problem and show numerical results for the quadrotor model.
We use the same system parameters as in \cite{madani2006backstepping}. We take the mass $m=2kg$, the gravity $g=9.81 m/s^2$, the moment of of inertia $J_x=J_y=J_z/2=1.2416 kg\cdot m^2$.
In the following experiments, we specify the domain of initial state as $\mathcal{S}_0=\{x,y\in [-10,10], z\in[5,100], v_x, v_y, v_z \in [-0.5, 0.5],\theta,\phi\in[-\pi/4, \pi/4],\psi\in[-\pi, \pi];\bm{w}=\bm{0}\}$. We always uniformly sample 100 samples of $\bm{x}_0$ from $\mathcal{S}_0$ to estimate the success rate and average computation time of the algorithm.
In Figure \ref{fig:uav_3d}, we present an example solution of the optimal landing problem that manipulates the quadrotor from the starting position $\bm{x}_0$ to the origin.
Figure \ref{subfig:path_3d} shows the optimal trajectory in the Earth-fixed coordinates obtained by the proposed method.
Figure \ref{subfig:p_3d} and Figure \ref{subfig:eta_3d} show the optimal position $\bm{p}$ and attitude $\bm{\eta}$, respectively.

In Section \ref{sec:baseline_solver} below, we first demonstrate that the baseline TPBVP solver for the free terminal time TPBVP~\eqref{eq:bvp_free_a}-\eqref{eq:bvp_free_c} hardly works. Then in Section \ref{sec:fix_initguess}, we show that the solution to the fixed terminal time TPBVP~\eqref{eq:bvp_landing}-\eqref{eq:bvp_bd_landing} can provide better initialization for the free terminal time problem. Finally, in Sections \ref{sec:space_marching} and \ref{sec:nn_predict}, we propose the space-marching technique and DNN-based prediction of the optimal terminal time to further improve the success rate and efficiency of the solver to the fixed time problem.
%

\subsection{Baseline method for the TPBVP}
\label{sec:baseline_solver}
Throughout the paper, we will use the classical \textit{bvp4c} method \citep{kierzenka2001bvp} as the TPBVP solver. 
It is a collocation method based on piecewise cubic polynomials. The solution at all collocation points and unknown parameters are solved from a system of algebraic equations by the quasi-Newton method. Then the residual defined in terms of the interpolant is used to estimate the error and refine the mesh. By the nature of the quasi-Newton method, the quality of the initial guess of the solution is critical for the performance.

In our TPBVP, without any prior knowledge, the simplest choice of the initial guess is to set $t_f$ to a reasonable scalar and $\bm{x}(t), \bm{\lambda}(t)$ to the constant zero, as summarized in Algorithm~\ref{alg:direct}. However, the TPBVP solver hardly converges with this choice. 
Table \ref{tab:direct} reports the success rate of Algorithm~\ref{alg:direct} with a few different initial guesses of $t_f$
We can see that, with zero initialization of the path, the success rate of solving the free terminal time TPBVP is always very low, regardless of the initial guess of $t_f$.

\begin{algorithm}
    \caption{Free terminal time problem with zero initialization}
  \begin{algorithmic}[1]
      \STATE \textbf{Input:} The initial state $\bm{x}_0$; the guess value of the optimal terminal time $\tilde{t}^*_f$.
    \STATE Solve the TPBVP corresponding to the free time problem with the zero as the initial guess of the path and $\tilde{t}^*_f$ as the initial guess of the terminal time. \\
    \STATE \textbf{Output:} The solution of the free time problem.
  \end{algorithmic} 
  \label{alg:direct}
\end{algorithm} 

\begin{table}[htpb]
\centering
\caption{Solving TPBVP corresponding to the free terminal time problems with zero initialization}
\medskip
\begin{tabular}{c|cccccc}
\toprule
    initial guess of $t_f$& 4 & 8 & 12 & 16 & 20 & 24\\ \midrule
success rate &3\% &4\% &0\%  &0\%  &1\% & 1\%\\ \bottomrule
\end{tabular}\label{tab:direct}
\end{table}

\subsection{Using solution to the fixed terminal time problem as initial guess}
\label{sec:fix_initguess}
Table \ref{tab:direct} suggests that such simple initial guess of $\bm{x}(t)$ and $\bm{\lambda}(t)$ leads to the TPBVP solver's poor performance when solving the free terminal time problem.
To address this issue, we notice that the solution of the free terminal time problem is also the solution of a corresponding fixed terminal time problem if the fixed terminal time $t_f$ equals the optimal terminal time $t^*_f$. In other words, if we have a reasonable guess of $t^*_f$, the solution of the fixed terminal time problem can provide us a good initial guess to the free terminal time problem. Moreover, the fixed terminal time problem is easier to solve with many efficient techniques, such as the marching method introduced in the following subsection.
Therefore, we can first guess a value of the optimal terminal time $\tilde{t}^*_f$ and solve the fixed terminal time problem with $t_f=\tilde{t}^*_f$. Then we use its solution as the initial guess to solve the free terminal time problem. This approach can be viewed as a warm start method for solving the optimization problem.
The corresponding algorithm is summarized in Algorithm \ref{alg:fix_as_init} (\textit{warm start with fixed terminal time solution}), in which the fixed time problem is solved using zero initialization.
\begin{algorithm}[!htb]
    \caption{Warm start with fixed terminal time solution}
  \begin{algorithmic}[1]
      \STATE \textbf{Input:} The initial state $\bm{x}_0$; the guess value of the optimal terminal time $\tilde{t}^*_f$.
    \STATE Solve the fixed time problem with $t_f=\tilde{t}^*_f$ using zero initialization. 
    \STATE Solve the free time problem by using the solution from step 2 as initial guess.    \STATE \textbf{Output:} The solution of the free terminal time problem.
  \end{algorithmic} 
  \label{alg:fix_as_init}
\end{algorithm} 

\begin{figure}[!htb]
\centering
\includegraphics[width=0.6\textwidth]{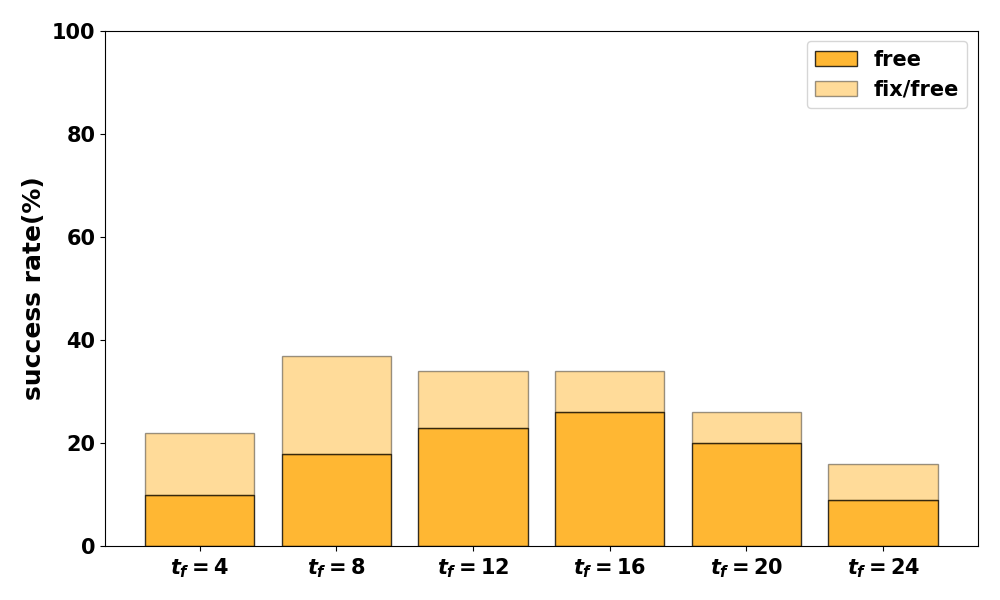}
\caption{Warm start with fixed terminal time solution: the label ``free" denotes the success rate for solving the free terminal time problems; the label ``fix/free" (with a lighter color) denotes the rate that the fixed terminal time problems are solved successfully but the corresponding free terminal time problem is not solved successfully.  The sum of the two parts is the success rate for solving the fixed terminal time problem.}
\label{fig:fix_as_warm}
\end{figure}

The numerical results of Algorithm \ref{alg:fix_as_init} with different choices of the initial guess $t_f=\tilde{t}^*_f$ are presented in Figure \ref{fig:fix_as_warm}.
Comparing Figure \ref{fig:fix_as_warm} with Table~\ref{tab:direct}, we can see that a warm start with the solution to the fixed terminal time problem significantly improves the success rate, although the rate is still not high enough for practical applications.

\subsection{Space-marching method for solving the fixed terminal time problem}
\label{sec:space_marching}
To further increase the success rate of solving the fixed terminal time TPBVP, we propose a homotopy method in the spirit of space-marching \citep{ascher1995numerical}. Space-marching is primarily developed to solve differential equations. Here we tailor a similar idea to help solve our optimal landing problems.
Its intuition is as follows.
Solving the fixed time problem is still difficult since the initial state $\bm{x}_0$ is far away from the terminal set $\mathcal{S}_{T}$.
We can solve a simpler fixed time problem whose initial state is closer to the terminal state and the corresponding solution is not far from that of the original fixed time problem. 
After the simpler fixed time problem is solved, we can use its solution as the initial guess to solve the original harder one. 
In other words, we seek another level of a warm start to help solve the fixed terminal time problem, and this process can be performed repeatedly.

To present this method in a systematic way,
we say $\bm{x}_{end}$ is a terminal state if there exists $\bm{u} \in \mathcal{U}$ such that for any $t_f \ge 0$,
the path
\begin{equation*}
    \bm{x}(t) \equiv \bm{x}_{end},\; \bm{u}(t) \equiv \bm{u},\; 0\leq t\leq t_f
\end{equation*}
is the optimal path for the fixed terminal time problem with $t_f$ as terminal time and $\bm{x}_0 = \bm{x}_{end}$. 
We always assume such a terminal state exists for the optimal landing problem. In this paper, we choose the origin (of the 12-dimensional state space) as the terminal state.
In order to solve the problem with a given initial state $\bm{x}_0$, we evenly select $K$ points in the line segment from $\bm{x}_{end}$ to $\bm{x}_0$, and denote them as $\{\bm{x}^1_0, \bm{x}^2_0, \cdots, \bm{x}^K_0\}$ according to their increasing distances to $\bm{x}_{end}$ ($\bm{x}^K_0 = \bm{x}_0$). 
We use $\bm{x}^{k}_{aug}$ to denote the solution $(\bm{x}(t),\bm{\lambda}(t)), (0\leq t\leq\tilde{t}_f^*)$ to the fixed time problem with the initial state $\bm{x}^{k}_{0}$ and terminal time $\tilde{t}_f^*$, $k=1,\dots,K$. We assume $\bm{x}^{0}_{aug}$ constant zero.
The space-marching starts with $k=1$.
In the $k$-\textit{th} step of marching, 
we solve the fixed time problem with the initial state $\bm{x}^k_0$ by using the solution $\bm{x}^{k-1}_{aug}$ obtained from the previous step as the initial guess.
The process repeats until $k=K$.
We call this algorithm \textit{warm start with fixed terminal time solution through space-marching}, and it is summarized in Algorithm~\ref{alg:space_marching}.

\begin{algorithm}[!htpb]
  \caption{Warm start with fixed terminal time solution through space-marching}
  \label{alg:SpaceMarch}
  \begin{algorithmic}[1]
      \STATE \textbf{Input:} The initial state $\bm{x}_0$; the guess value of the optimal terminal time $\tilde{t}^*_f$; the the number of marching steps $K$.
    \STATE Evenly select $K$ points in the line segment from $\bm{x}_{end}$ to $\bm{x}_0$, and denotes them as $\{\bm{x}_0^1, \bm{x}_0^2,\cdots,\bm{x}_0^K\}$.
    \STATE Initialize $\bm{x}^{0}_{aug}$ with constant zero.
    \FOR{$k=1,2,\cdots,K$}
    \STATE Solve the fixed terminal time problem with initial state $\bm{x}_0^k$ and terminal time $\tilde{t}^*_f$ by using $\bm{x}^{k-1}_{aug}$ as the initial guess. Denote the solution as $\bm{x}^{k}_{aug}$.
    \ENDFOR \\
    Solving the free time problem with $\bm{x}^{K}_{aug}$ as the initial guess of the path and $\tilde{t}^*_f$ as the initial guess of the terminal time.
    \STATE \textbf{Output:} The solution of the free terminal time problem with initial state $\bm{x}_0$.
  \end{algorithmic}
  \label{alg:space_marching}
\end{algorithm}

\begin{figure}[h!]
\centering
\includegraphics[width=0.75\textwidth]{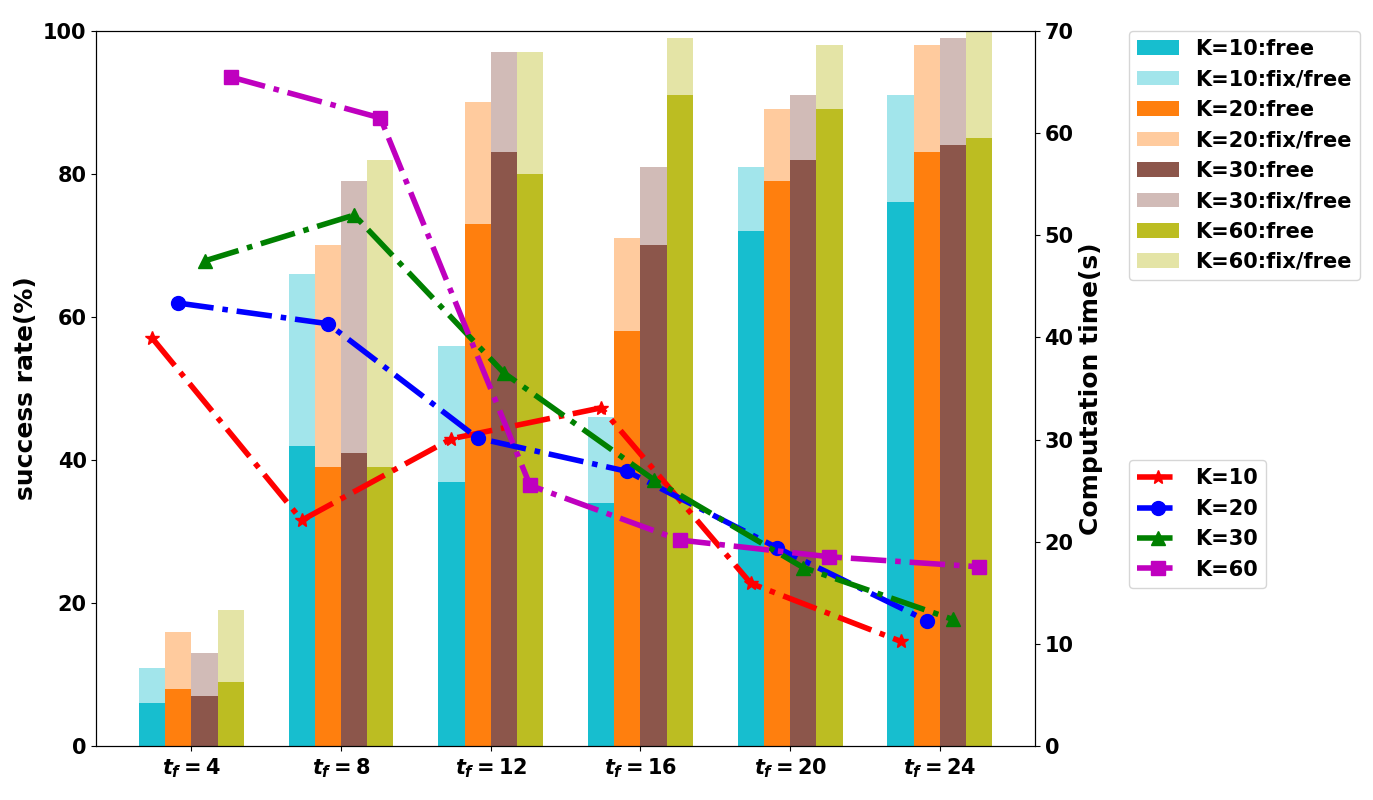}
\caption{Warm start with fixed terminal time solution through space-marching: the histograms denote the success rates (left $y$-axis), whose labels share the same meaning as those in Figure~\ref{fig:fix_as_warm}; the dash-dotted lines denote the average computation time (in seconds) per path (right $y$-axis).}
\label{fig:fix_as_warm_march}
\end{figure}

\begin{figure}[h!]
\centering
\includegraphics[width=0.55\textwidth]{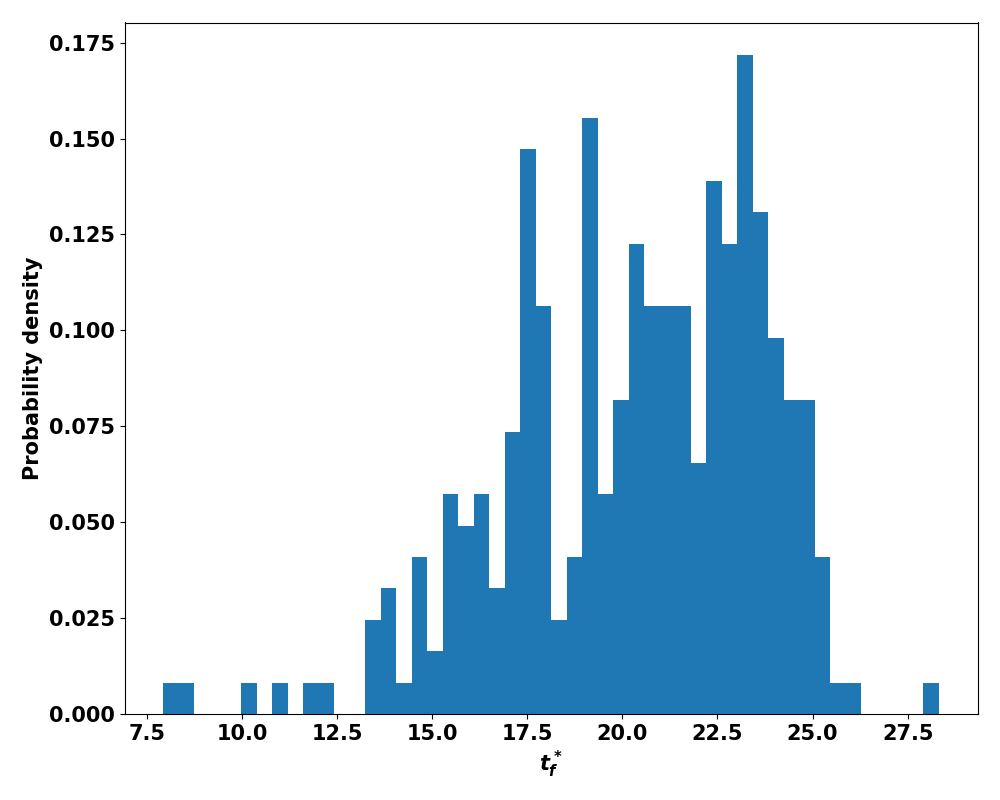}
\caption{The distribution of optimal terminal time $t^*_f$ obtained by solving 300 optimal landing problems with randomly selected initial states.} 
\label{fig:tf_dist}
\end{figure}

Figure~\ref{fig:fix_as_warm_march} shows the success rate and computation time of Algorithm~\ref{alg:space_marching} with different choices of the initial guess $t_f=\tilde{t}^*_f$ and marching steps $K$.
We can see that most of the fixed terminal time problems can be solved with a high success rate if the guessed terminal time $\tilde{t}_f^*$ is not too small and the number of marching steps $K$ is large enough (greater than 20). However, the free terminal time problem might not be solved successfully if the guessed terminal time is not close enough to the optimal terminal time.
To see this more clearly, we use Algorithm~\ref{fig:fix_as_warm_march} with $K=60, \tilde{t}^*_f=24$ to collect optimal time $t_f^*$ associated with 300 randomly sampled initial states $\bm{x}_0$.
We plot the distribution of $t_f^*$ in Figure~\ref{fig:tf_dist}, from which we can see that the values of $t_f^*$ are distributed over a wide range. Hence, using a constant initial guess $\tilde{t}_f^*$ cannot achieve the best performance. To improve the success rate of the free terminal time problem, we need a more accurate prediction of the optimal terminal time. We will investigate this issue in the next subsection.

\subsection{Predicting the optimal terminal time}
\label{sec:nn_predict}
In this section, we consider empowering Algorithm~\ref{alg:space_marching} by predicting the optimal terminal time as a function of the initial state $\bm{x}_0$ through a linear model or a neural network.
To do so, we need to prepare a dataset for supervised learning. As described in the previous section, we first randomly select 300 initial states $\bm{x}_0$ and then use Algorithm~\ref{fig:fix_as_warm_march} with $K=60, \tilde{t}^*_f=24$ to collect 300 optimal landing paths (the specified terminal time 24 is not necessarily the optimal terminal time).
We select 100 states (uniformly in time) on each optimal landing path and store the corresponding optimal landing time to obtain the training data. We have 30000 pairs of starting positions and optimal ending times for training in total.
We then use this dataset to optimize a linear model and a neural network model (3 three hidden layers and 64 neurons in each layer) based on the objective being the squared difference between the predicted $\tilde{t}_f^*$ and the truth optimal terminal time. The Adam optimizer \citep{kingma2014adam} is used to train the neural network model with 500 epochs, batch size 128, and learning rate 0.002. Afterwards, when we need to solve a free terminal time problem with a new initial state $\bm{x}_0$, we first use the linear model or the neural network to predict the optimal terminal time $\tilde{t}_f^*$ associated with $\bm{x}_0$ and then use Algorithm~\ref{alg:space_marching} to solve the problem.

\begin{figure}[h!]
\centering
\includegraphics[width=0.65\textwidth]{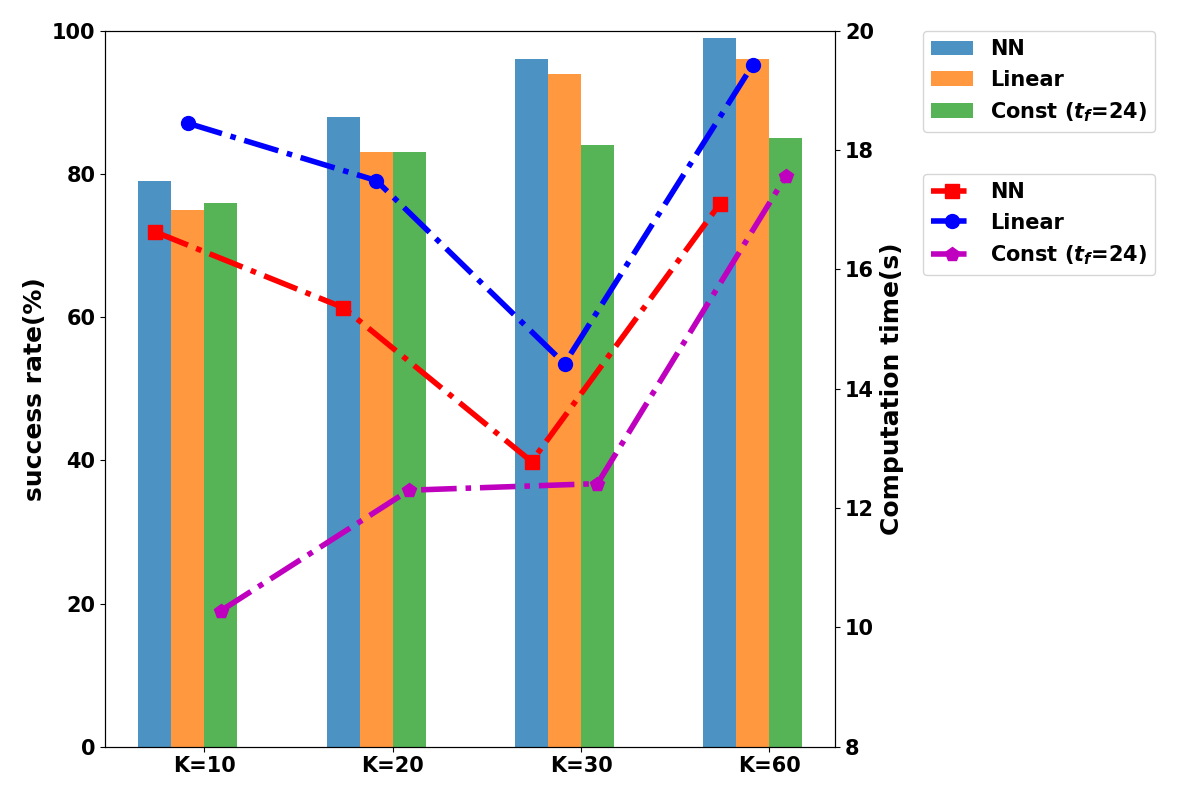}
\caption{Comparison between the constant, linear model, and neural network model when they are used to predict $t^*_f$ as the initial guess in the TPBVP solver.  The histograms denote the success rates (left $y$-axis) of solving the free terminal time problem. Dash-dotted lines denote the average computation time (in seconds) per path (right $y$-axis).}
\label{fig:NN_vs_linear}
\end{figure}

The success rates of using a constant $(\tilde{t}_f^*=24)$, linear model, or a neural network to predict the optimal terminal time with the different space-marching steps $K$ are presented in Figure~\ref{fig:NN_vs_linear}. Comparing these results, together with those using other guessed constant optimal terminal time in Figure \ref{fig:fix_as_warm_march}, we can see that both the linear model and neural network model achieve much higher success rates. Using a neural network attains higher success rates and takes less computation time because it can predict the optimal terminal time more accurately. With the help of neural networks and the space-marching with $K = 60$, we achieve a $99\%$ successful rate, and the average computational time is about 17 seconds, which is the best performance among the methods considered in this paper.

\section{Suboptimality of Model Predictive Control}
\label{sec:MPC}
In this section, we study the performance of the model predictive control (MPC) method \citep{camacho2013model}, a widely used algorithm in optimal control, on the optimal landing problem.
At each discretized time step, the algorithm computes a cost-minimizing control strategy starting from the current state for a relatively short time horizon (called the prediction horizon). 
The obtained control strategy is implemented for a very small time duration, then the states are updated and the calculations are repeated.
There have been various studies using the MPC method to control quadrotors \citep{bangura2014real,ru2017nonlinear,eren2017model,romero2021model}.
For application to the landing problem \citep{eren2017model}, a landing trajectory is required as the input. MPC can then use that to define the running cost in the prediction horizon in order to track that trajectory to land on the target position.
%


\begin{figure}[!htb]
	\centering
	\begin{subfigure}[b]{.7\textwidth}
	\includegraphics[width=\textwidth]{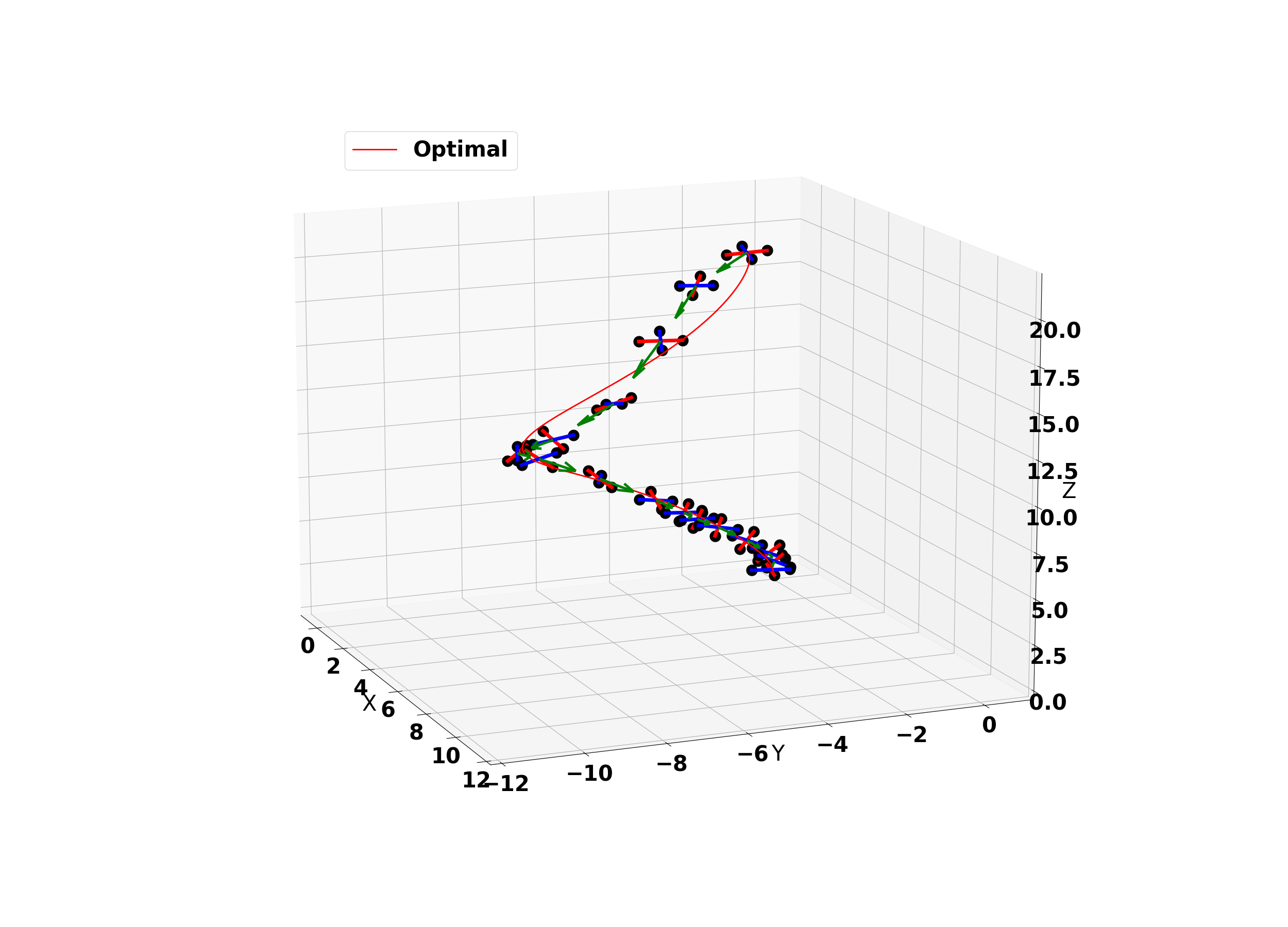}
	\vspace{-3.5em}
	\caption{}
	\label{subfig:path_mpc}
	\end{subfigure}
	\begin{subfigure}[b]{.85\textwidth}
	\includegraphics[width=\textwidth]{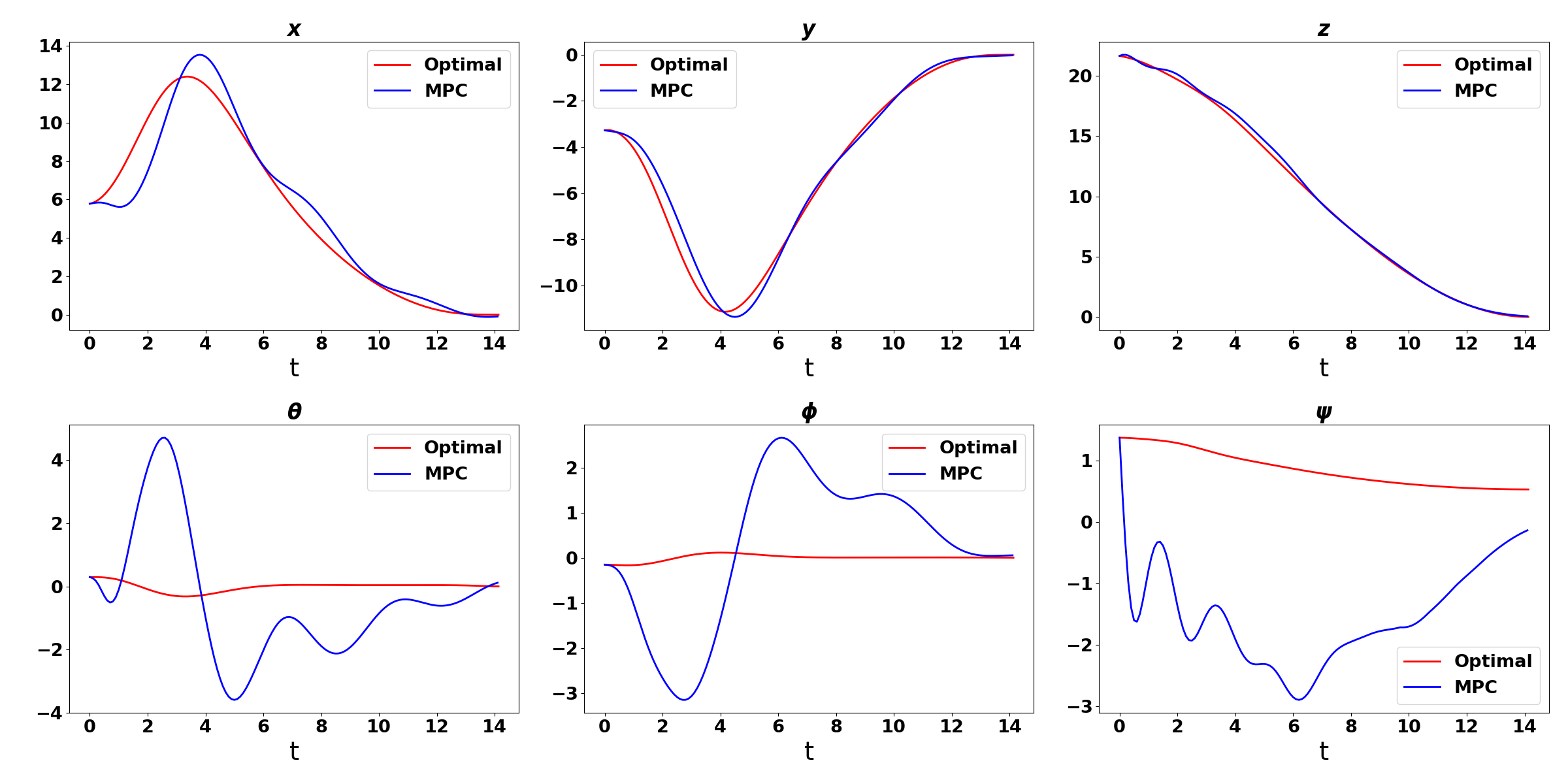}
	\caption{}
	\label{subfig:path_compare}
	\end{subfigure}
	\caption{(a): The trajectory obtained by the MPC method (cost=32.31) following the optimal trajectory obtained by the proposed method (cost=18.63); (b): The optimal position $(x, y, z)$ and attitude $(\phi, \theta, \psi)$ (red) vs. those obtained by the MPC method (blue).}
	\label{fig:mpc}
\end{figure}

We pick an initial state and compute its optimal trajectory with our machine learning enhanced algorithm. The overall trajectory and trajectories of position and attitude are shown in Figure \ref{subfig:path_mpc} and \ref{subfig:path_compare} with red color. 
The corresponding optimal terminal time $t_f^*$ is 14.58(s), and the optimal cost obtained by the proposed method is 18.63.
We then use the Matlab Model Predictive Control Toolbox \citep{bemporad2004model} to track the optimal trajectory with the MPC method.
We set the prediction horizon as $1.8$(s) with a time difference $\Delta t = 0.1$(s) to strike a reasonable balance between computation accuracy and efficiency. We also modify the cost parameter in the prediction horizon to achieve the best tracking performance.
The blue curves in Figure \ref{subfig:path_compare} show the trajectories of the position and attitude obtained by the MPC method. The corresponding total cost according to \eqref{eq:cost} is 32.31. From this result, we can clearly see that the MPC method can only find controls with suboptimality even provided the optimal trajectory. It is determined by the nature of the short horizon approximation scheme in the MPC method.

\section{Conclusion and discussion}
\label{sec:summary}
This paper proposes a machine learning enhanced algorithm for solving the optimal landing problem.
The proposed algorithm is based on providing a good initialization for the corresponding TPBVP derived from PMP.
It has two main components: predicting the terminal time using DNN and solving the resulting fixed time problem using a space-marching method.
Through a series of experiments, the effectiveness of the proposed method has been verified. 

\begin{figure}[!htb]
	\centering
	\begin{subfigure}[b]{.28\textwidth}
	\includegraphics[width=\textwidth]{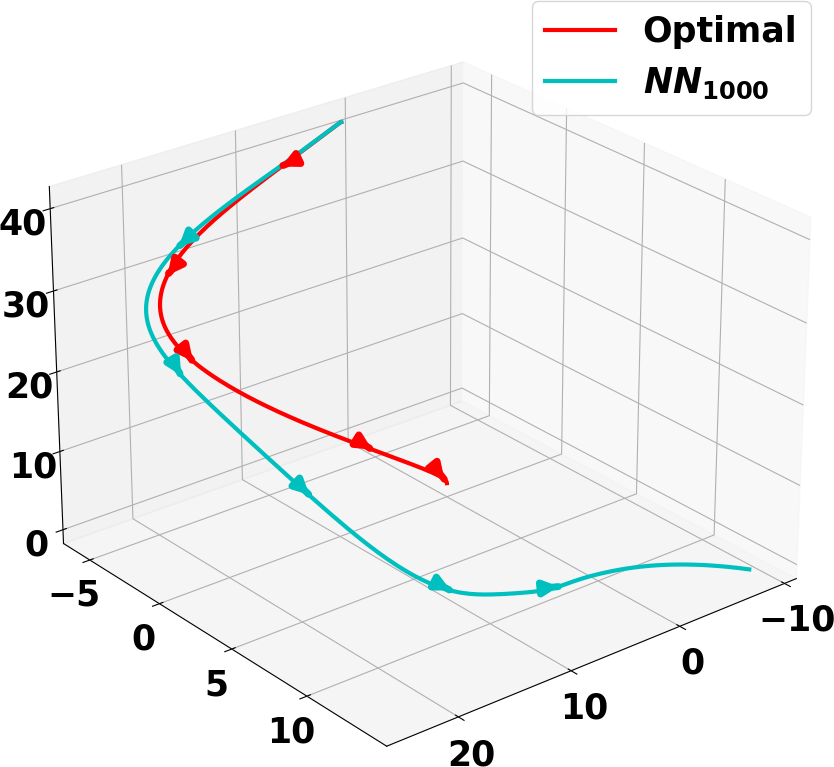}
	\caption{$NN_{1000}$}
	\label{subfig:nn_1000}
	\end{subfigure}
	~~~~
	\begin{subfigure}[b]{.28\textwidth}
	\includegraphics[width=\textwidth]{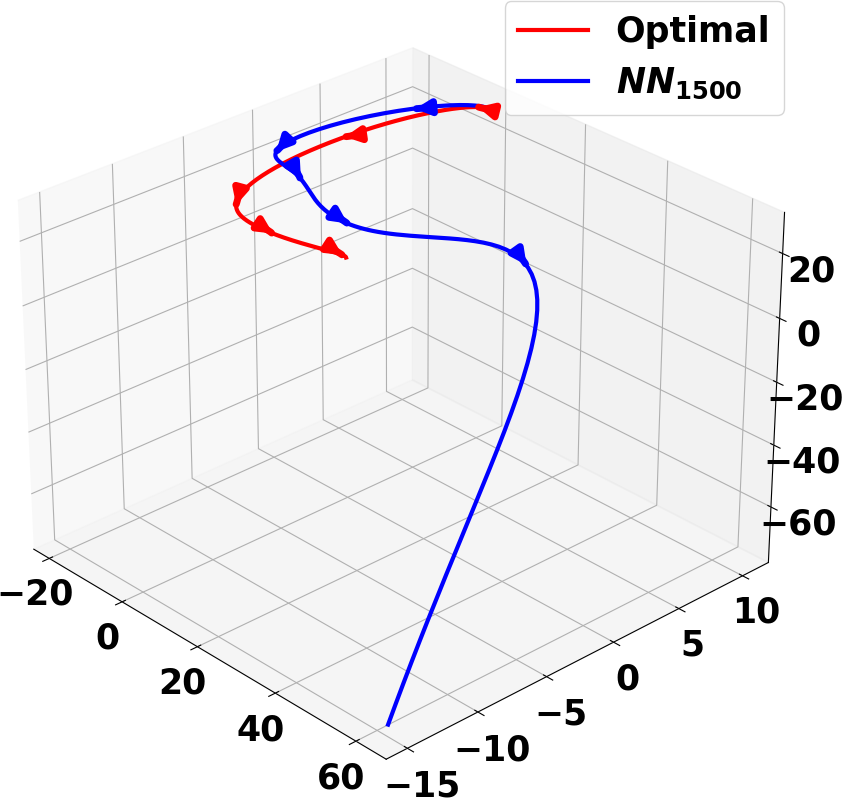}
	\caption{$NN_{1500}$}
	\label{subfig:nn_1500}
	\end{subfigure}
	~~~~
	\begin{subfigure}[b]{.28\textwidth}
	\includegraphics[width=\textwidth]{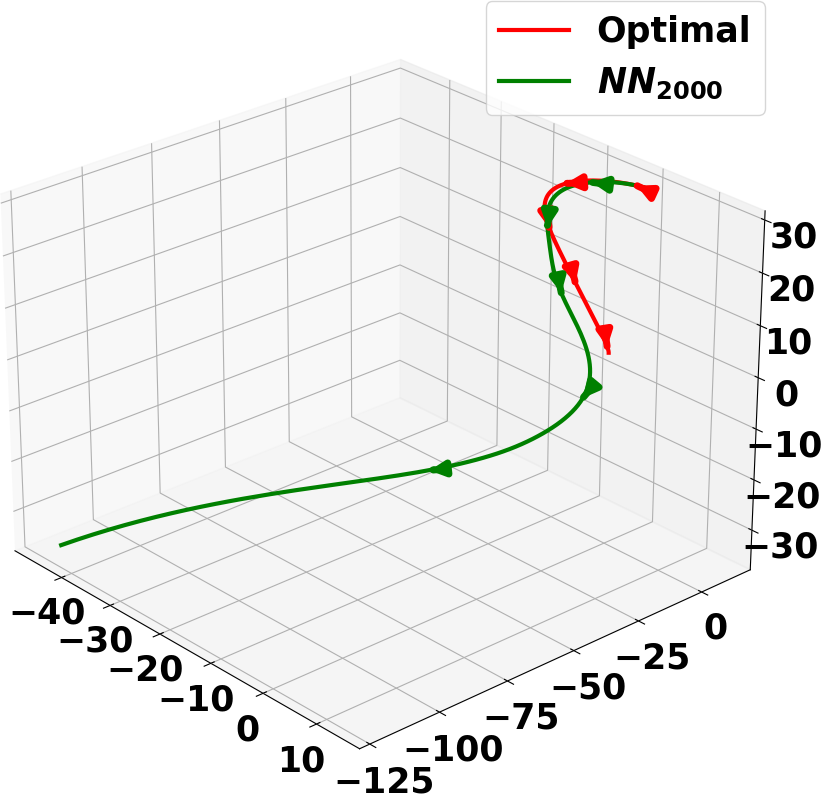}
	\caption{$NN_{2000}$}
	\label{subfig:nn_2000}
	\end{subfigure}
	\caption{Trajectories on the time interval $[0, t_f^*]$, starting from the same initial state and following controls provided by neural networks $NN_{1000}$, $NN_{1500}$, and $NN_{2000}$, respectively. The neural network-based feedback control performs reasonably well at the beginning and then deteriorates.}
	\label{fig:nn_v}
\end{figure}

So far we have mainly used neural networks to predict the optimal terminal time.
An obvious alternative is to directly learn the mapping from the state to control or value function using DNN.
Our effort along this line has not been successful so far.
Here we present such examples.
We consider the same problem in Section \ref{sec:DNN}. Similarly as above, we randomly sample 1000, 1500, and 2000 samples of $\bm{x}_0$ from $\mathcal{S}_0$ and use them to generate 1000, 1500, and 2000 optimal trajectories, respectively. We then collect 20 (uniformly in time) states on each path and the corresponding value, resulting in three training sets with 20000, 30000, and 40000 pairs of data, respectively. Finally, we use these training sets to train three different neural networks, denoted by $NN_{1000}$, $NN_{1500}$, and $NN_{2000}$, respectively. According to the PMP, the gradients of the approximate value function should provide the optimal feedback control~\citep{nakamura2021adaptive}. Figures \ref{subfig:nn_1000}, \ref{subfig:nn_1500} and \ref{subfig:nn_2000} show trajectories on the time interval $[0, t_f^*]$, starting from the same initial state and following such control strategies provided by three neural networks.
Compared to the optimal trajectory obtained by our algorithm, the neural network-based feedback control performs reasonably in the early stages of the landing but gradually deteriorates.
Even when we use these trajectories as the initial guess for the TPBVP solver, the solution still does not converge.
This phenomenon is perhaps due to the high dimensionality and strong nonlinearity of the problem, as well as the wide range of scales spanned by
  the training data.
It is an important future direction to explore techniques that can leverage neural networks to directly provide feedback control robustly for challenging problems like the optimal landing problem considered in this paper. 

\bibliographystyle{abbrvnat}
\bibliography{uav_ref}

\begin{thebibliography}{39}
\providecommand{\natexlab}[1]{#1}
\providecommand{\url}[1]{\texttt{#1}}
\expandafter\ifx\csname urlstyle\endcsname\relax
  \providecommand{\doi}[1]{doi: #1}\else
  \providecommand{\doi}{doi: \begingroup \urlstyle{rm}\Url}\fi

\bibitem[Ackermann(2020)]{ackermann2020ai}
E.~Ackermann.
\newblock {AI}-powered drone learns extreme acrobatics.
\newblock \emph{IEEE Spectrum}, 2020.

\bibitem[Ascher et~al.(1995)Ascher, Mattheij, and Russell]{ascher1995numerical}
U.~M. Ascher, R.~M. Mattheij, and R.~D. Russell.
\newblock \emph{Numerical solution of boundary value problems for ordinary
  differential equations}.
\newblock SIAM, 1995.

\bibitem[Assellaou et~al.(2016)Assellaou, Bokanowski, Desilles, and
  Zidani]{assellaou2016hamilton}
M.~Assellaou, O.~Bokanowski, A.~Desilles, and H.~Zidani.
\newblock A {H}amilton-{J}acobi-{B}ellman approach for the optimal control of
  an abort landing problem.
\newblock In \emph{2016 IEEE 55th Conference on Decision and Control (CDC)},
  pages 3630--3635. IEEE, 2016.

\bibitem[Bangura and Mahony(2014)]{bangura2014real}
M.~Bangura and R.~Mahony.
\newblock Real-time model predictive control for quadrotors.
\newblock \emph{IFAC Proceedings Volumes}, 47\penalty0 (3):\penalty0
  11773--11780, 2014.

\bibitem[Bemporad et~al.(2004)Bemporad, Morari, and Ricker]{bemporad2004model}
A.~Bemporad, M.~Morari, and N.~L. Ricker.
\newblock Model predictive control toolbox.
\newblock \emph{User's Guide, Version}, 2, 2004.

\bibitem[Beppu et~al.(2021)Beppu, Maruta, and Fujimoto]{beppu2021value}
H.~Beppu, I.~Maruta, and K.~Fujimoto.
\newblock Value iteration with deep neural networks for optimal control of
  input-affine nonlinear systems.
\newblock \emph{SICE Journal of Control, Measurement, and System Integration},
  14\penalty0 (1):\penalty0 140--149, 2021.

\bibitem[Bertolazzi et~al.(2005)Bertolazzi, Biral, and
  Da~Lio]{bertolazzi2005symbolic}
E.~Bertolazzi, F.~Biral, and M.~Da~Lio.
\newblock Symbolic--numeric indirect method for solving optimal control
  problems for large multibody systems.
\newblock \emph{Multibody System Dynamics}, 13\penalty0 (2):\penalty0 233--252,
  2005.

\bibitem[Bock and Plitt(1984)]{bock1984multiple}
H.~G. Bock and K.-J. Plitt.
\newblock A multiple shooting algorithm for direct solution of optimal control
  problems.
\newblock \emph{IFAC Proceedings Volumes}, 17\penalty0 (2):\penalty0
  1603--1608, 1984.

\bibitem[Bouktir et~al.(2008)Bouktir, Haddad, and
  Chettibi]{bouktir2008trajectory}
Y.~Bouktir, M.~Haddad, and T.~Chettibi.
\newblock Trajectory planning for a quadrotor helicopter.
\newblock In \emph{2008 16th mediterranean conference on control and
  automation}, pages 1258--1263. Ieee, 2008.

\bibitem[Camacho and Alba(2013)]{camacho2013model}
E.~F. Camacho and C.~B. Alba.
\newblock \emph{Model predictive control}.
\newblock Springer science \& business media, 2013.

\bibitem[Chen et~al.(2018)Chen, Shi, and Zhang]{chen2018optimal}
Y.~Chen, Y.~Shi, and B.~Zhang.
\newblock Optimal control via neural networks: A convex approach.
\newblock \emph{arXiv preprint arXiv:1805.11835}, 2018.

\bibitem[Eren et~al.(2017)Eren, Prach, Ko{\c{c}}er, Rakovi{\'c}, Kayacan, and
  A{\c{c}}{\i}kme{\c{s}}e]{eren2017model}
U.~Eren, A.~Prach, B.~B. Ko{\c{c}}er, S.~V. Rakovi{\'c}, E.~Kayacan, and
  B.~A{\c{c}}{\i}kme{\c{s}}e.
\newblock Model predictive control in aerospace systems: Current state and
  opportunities.
\newblock \emph{Journal of Guidance, Control, and Dynamics}, 40\penalty0
  (7):\penalty0 1541--1566, 2017.

\bibitem[Foehn et~al.(2021)Foehn, Romero, and Scaramuzza]{foehn2021time}
P.~Foehn, A.~Romero, and D.~Scaramuzza.
\newblock Time-optimal planning for quadrotor waypoint flight.
\newblock \emph{Science Robotics}, 6\penalty0 (56):\penalty0 eabh1221, 2021.

\bibitem[Geisert and Mansard(2016)]{geisert2016trajectory}
M.~Geisert and N.~Mansard.
\newblock Trajectory generation for quadrotor based systems using numerical
  optimal control.
\newblock In \emph{2016 IEEE international conference on robotics and
  automation (ICRA)}, pages 2958--2964. IEEE, 2016.

\bibitem[Guo et~al.(2011)Guo, Hawkins, and Wie]{guo2011optimal}
Y.~Guo, M.~Hawkins, and B.~Wie.
\newblock Optimal feedback guidance algorithms for planetary landing and
  asteroid intercept.
\newblock In \emph{AAS/AIAA astrodynamics specialist conference}, pages
  2011--588. AAS, 2011.

\bibitem[Han and E(2016)]{han2016deep}
J.~Han and W.~E.
\newblock Deep learning approximation for stochastic control problems.
\newblock \emph{arXiv preprint arXiv:1611.07422}, 2016.

\bibitem[Hartl et~al.(1995)Hartl, Sethi, and Vickson]{hartl1995survey}
R.~F. Hartl, S.~P. Sethi, and R.~G. Vickson.
\newblock A survey of the maximum principles for optimal control problems with
  state constraints.
\newblock \emph{SIAM review}, 37\penalty0 (2):\penalty0 181--218, 1995.

\bibitem[Hu and Mishra(2017)]{hu2017time}
B.~Hu and S.~Mishra.
\newblock A time-optimal trajectory generation algorithm for quadrotor landing
  onto a moving platform.
\newblock In \emph{2017 American Control Conference (ACC)}, pages 4183--4188.
  IEEE, 2017.

\bibitem[Hu et~al.(2015)Hu, Lu, and Mishra]{hu2015fast}
B.~Hu, L.~Lu, and S.~Mishra.
\newblock Fast, safe and precise landing of a quadrotor on an oscillating
  platform.
\newblock In \emph{2015 American Control Conference (ACC)}, pages 3836--3841.
  IEEE, 2015.

\bibitem[Kierzenka and Shampine(2001)]{kierzenka2001bvp}
J.~Kierzenka and L.~F. Shampine.
\newblock A {BVP} solver based on residual control and the {Maltab} {PSE}.
\newblock \emph{ACM Transactions on Mathematical Software (TOMS)}, 27\penalty0
  (3):\penalty0 299--316, 2001.

\bibitem[Kingma and Ba(2015)]{kingma2014adam}
D.~P. Kingma and J.~Ba.
\newblock Adam: a method for stochastic optimization.
\newblock In \emph{Proceedings of the International Conference on Learning
  Representations}, 2015.

\bibitem[Liu(2019)]{liu2019fuel}
X.~Liu.
\newblock Fuel-optimal rocket landing with aerodynamic controls.
\newblock \emph{Journal of Guidance, Control, and Dynamics}, 42\penalty0
  (1):\penalty0 65--77, 2019.

\bibitem[Liu and Duan(2006)]{liu2006nonlinear}
X.~Liu and G.~Duan.
\newblock Nonlinear optimal control for the soft landing of lunar lander.
\newblock In \emph{2006 1st International Symposium on Systems and Control in
  Aerospace and Astronautics}, pages 6--pp. IEEE, 2006.

\bibitem[Madani and Benallegue(2006)]{madani2006backstepping}
T.~Madani and A.~Benallegue.
\newblock Backstepping control for a quadrotor helicopter.
\newblock In \emph{2006 IEEE/RSJ International Conference on Intelligent Robots
  and Systems}, pages 3255--3260. IEEE, 2006.

\bibitem[Nakamura-Zimmerer et~al.(2021{\natexlab{a}})Nakamura-Zimmerer, Gong,
  and Kang]{nakamura2021adaptive}
T.~Nakamura-Zimmerer, Q.~Gong, and W.~Kang.
\newblock Adaptive deep learning for high-dimensional
  {H}amilton--{J}acobi--{B}ellman equations.
\newblock \emph{SIAM Journal on Scientific Computing}, 43\penalty0
  (2):\penalty0 A1221--A1247, 2021{\natexlab{a}}.

\bibitem[Nakamura-Zimmerer et~al.(2021{\natexlab{b}})Nakamura-Zimmerer, Gong,
  and Kang]{nakamura2021neural}
T.~Nakamura-Zimmerer, Q.~Gong, and W.~Kang.
\newblock Neural network optimal feedback control with enhanced closed loop
  stability.
\newblock \emph{arXiv preprint arXiv:2109.07466}, 2021{\natexlab{b}}.

\bibitem[Romero et~al.(2021)Romero, Sun, Foehn, and
  Scaramuzza]{romero2021model}
A.~Romero, S.~Sun, P.~Foehn, and D.~Scaramuzza.
\newblock Model predictive contouring control for near-time-optimal quadrotor
  flight.
\newblock \emph{arXiv preprint arXiv:2108.13205}, 2021.

\bibitem[Ross and Fahroo(2002)]{ross2002direct}
I.~M. Ross and F.~Fahroo.
\newblock A direct method for solving nonsmooth optimal control problems.
\newblock \emph{IFAC Proceedings Volumes}, 35\penalty0 (1):\penalty0 479--484,
  2002.

\bibitem[Ru and Subbarao(2017)]{ru2017nonlinear}
P.~Ru and K.~Subbarao.
\newblock Nonlinear model predictive control for unmanned aerial vehicles.
\newblock \emph{Aerospace}, 4\penalty0 (2):\penalty0 31, 2017.

\bibitem[S{\'a}nchez-S{\'a}nchez and Izzo(2018)]{sanchez2018real}
C.~S{\'a}nchez-S{\'a}nchez and D.~Izzo.
\newblock Real-time optimal control via deep neural networks: study on landing
  problems.
\newblock \emph{Journal of Guidance, Control, and Dynamics}, 41\penalty0
  (5):\penalty0 1122--1135, 2018.

\bibitem[Shi et~al.(2019)Shi, Shi, O’Connell, Yu, Azizzadenesheli,
  Anandkumar, Yue, and Chung]{Shi2019Neural}
G.~Shi, X.~Shi, M.~O’Connell, R.~Yu, K.~Azizzadenesheli, A.~Anandkumar,
  Y.~Yue, and S.-J. Chung.
\newblock Neural lander: Stable drone landing control using learned dynamics.
\newblock In \emph{2019 International Conference on Robotics and Automation
  (ICRA)}, pages 9784--9790. IEEE, 2019.

\bibitem[Tang and Hauser(2019)]{tang2019data}
G.~Tang and K.~Hauser.
\newblock A data-driven indirect method for nonlinear optimal control.
\newblock \emph{Astrodynamics}, 3\penalty0 (4):\penalty0 345--359, 2019.

\bibitem[Tang et~al.(2018)Tang, Sun, and Hauser]{tang2018learning}
G.~Tang, W.~Sun, and K.~Hauser.
\newblock Learning trajectories for real-time optimal control of quadrotors.
\newblock In \emph{2018 IEEE/RSJ International Conference on Intelligent Robots
  and Systems (IROS)}, pages 3620--3625. IEEE, 2018.

\bibitem[Tsiotras et~al.(2011)Tsiotras, Bakolas, and Zhao]{tsiotras2011initial}
P.~Tsiotras, E.~Bakolas, and Y.~Zhao.
\newblock Initial guess generation for aircraft landing trajectory
  optimization.
\newblock In \emph{AIAA Guidance, Navigation, and Control Conference}, page
  6689, 2011.

\bibitem[Vanderbei(1999)]{vanderbei1999loqo}
R.~J. Vanderbei.
\newblock {LOQO}: An interior point code for quadratic programming.
\newblock \emph{Optimization methods and software}, 11\penalty0 (1-4):\penalty0
  451--484, 1999.

\bibitem[Verbeke and Schutter(2018)]{verbeke2018experimental}
J.~Verbeke and J.~D. Schutter.
\newblock Experimental maneuverability and agility quantification for rotary
  unmanned aerial vehicle.
\newblock \emph{International Journal of Micro Air Vehicles}, 10\penalty0
  (1):\penalty0 3--11, 2018.

\bibitem[Wang(2009)]{wang2009solving}
X.~Wang.
\newblock Solving optimal control problems with {MATLAB}: Indirect methods.
\newblock Technical report, 2009.

\bibitem[You et~al.(2020)You, Wan, Dai, Lu, and Rea]{you2020learning}
S.~You, C.~Wan, R.~Dai, P.~Lu, and J.~R. Rea.
\newblock Learning-based optimal control for planetary entry, powered descent
  and landing guidance.
\newblock In \emph{AIAA Scitech 2020 Forum}, page 0849, 2020.

\bibitem[Zhu et~al.(2019)Zhu, Ma, and Wang]{zhu2019deep}
L.~Zhu, J.~Ma, and S.~Wang.
\newblock Deep neural networks based real-time optimal control for lunar
  landing.
\newblock In \emph{IOP Conference Series: Materials Science and Engineering},
  volume 608, page 012045. IOP Publishing, 2019.

\end{thebibliography}

\end{document}